\documentclass[10pt]{amsart}
\usepackage{amssymb,amsmath}

\usepackage{float}
\usepackage{subfig,epsfig}
\usepackage{graphicx}
\usepackage{blindtext}
\usepackage{flafter}
\usepackage[export]{adjustbox}
\usepackage{caption}
\usepackage[backend=biber,
style=numeric
]{biblatex}
\usepackage{tensor}
\usepackage{cancel}
\newtheorem{proposition}{Proposition}[section]
\newtheorem{conjecture}{Conjecture}[section]
\newcommand{\bra}{\langle}
\newcommand{\ket}{\rangle}
\newcommand{\diff}[1]{\, d#1}

\newcommand{\hess}{\nabla^2}
\newcommand{\Rc}{\mbox{Rc}}

\usepackage{hyperref}
\hypersetup{
  colorlinks   = true,    
  urlcolor     = blue,    
  linkcolor    = blue,    
  citecolor    = red      
}

\newtheorem{theorem}{Theorem}[section]
\newtheorem{lemma}[theorem]{Lemma}

\theoremstyle{definition}

\theoremstyle{remark}

\numberwithin{equation}{section}

\begin{document}
	\title
	{Vanishing Bach-Like Tensors on Complete Gradient Shrinking Ricci Solitons}
	
	\author{James Siene}
	\address[James Siene]{Department of Mathematics\\ Lehigh University\\
		Bethlehem, PA 18015} 
	\email{jts614@lehigh.edu; }
	

\begin{abstract}
	The Bach tensor is classically defined in dimension 4, and work from J. Bergman \cite{bergman:2004} and others shows that $B = \frac{1}{2}U + \frac{1}{6}V$ where $U$ and $V$ are more basic 2-tensors, which are symmetric, divergence-free, algebraically independent, and quadratic in the Riemann tensor. In this paper, we extend H.-D. Cao and Q. Chen's results \cite{caochen:2013} for Bach-flat gradient shrinking Ricci solitons to solitons with $\mathfrak{B} = \alpha U + \beta V = 0$.
\end{abstract}

\maketitle
\date{}

\section{Notation and Conventions}
We begin by fixing some notations and commonly used tensors. We take the following conventions, assuming $\nabla$ is the Levi-Civita connection and $\Gamma_{ijk}$ are the Christoffel symbols:
\begin{align*}
    \nabla_i X^j &= \frac{\partial X^j}{\partial x_i} + \Gamma_{ip}^j X^p \\
    \nabla_i \omega_j &= \frac{\partial \omega_j}{\partial x_i} - \Gamma_{ij}^p \omega_p \\
    \nabla_i h_{jk} &= \frac{\partial h_{jk}}{\partial x_i} - \Gamma_{ij}^p h_{pk}- \Gamma_{ik}^p h_{jp}
\end{align*}
And the Ricci formula for commuting derivatives:
\begin{equation}\label{eq.commutingderivatives}
    \nabla_i \nabla_j X_k - \nabla_j \nabla_i X_k = -R_{ijk}^{\,\,\,l} X_l
\end{equation}
The Weyl tensor (or Weyl curvature) is given by:
\begin{equation}\label{eq.weyltensor}
    \begin{split}
        W_{ijkl} = &R_{ijkl} - \frac{1}{n-2}\left(g_{ik}R_{jl} - g_{il} R_{jk} - g_{jk}R_{il} + g_{jl}R_{ik}\right) \\
        &+ \frac{R}{(n-1)(n-2)}\left(g_{ik}g_{jl} - g_{il}g_{jk}\right)
    \end{split}
\end{equation}
The Cotton tensor plays a special role in the geometry of 3-dimensional manifolds. It is also useful for our purposes in dimension 4, so we list it here:
\begin{equation}\label{eq.cottontensor}
    C_{ijk} = \nabla_i R_{jk} - \nabla_j R_{ik} - \frac{1}{2(n-1)}\left(g_{jk}\nabla_i R - g_{ik} \nabla_j R\right)
\end{equation}
and it's a well-known fact that the Cotton tensor is the divergence of the Weyl tensor in dimensions $n > 1$:
\begin{equation}
    C_{ijk} = -\frac{n-2}{n-3}\nabla_l W_{ijkl}
\end{equation}
Moreover, the $D$-tensor is just the Cotton tensor of the conformal metric $\hat{g} = e^{\frac{2}{2-n}f} g$, so that:
\begin{equation*}
    D_{ijk} = C_{ijk} + W_{ijkl} \nabla_l f
\end{equation*}

\section{Motivation for Studying U and V on Solitons}
The Bach tensor was introduced by R. Bach \cite{bach:1921} for dimension $n=4$ in the 1920s to study conformal relativity and is given by:
\begin{equation}\label{eq.bachtensor}
    B_{ij} = \nabla^k\nabla^l W_{ikjl} + \frac{1}{2}R_{kl}\tensor{W}{_i^k_j^l}
\end{equation}
One could also see, and it is well known, that a metric is Bach-flat if and only if it is a critical point of the $S_W$-functional in dimension 4.
The Bach tensor is conformally invariant, divergence-free, and trace-free, as can be seen with direct calculation. Until 1968, it was the only known tensor with these properties that is also algebraically independent of the Weyl tensor. It is also quadratic in the Riemann tensor: this means that it can be expressed as products and contractions of exactly two copies of $R_m, Rc, $ and $R$, or second derivatives of one of these tensors. For example, the following expression is quadratic in the Riemann tensor:
\begin{equation*}
    R_{ij} R^{ij} + \frac{1}{2} R \, R_{ij} - \nabla_i \nabla_j R + \Delta R_{ij}
\end{equation*}
H.-D. Cao and Q. Chen used the following definition of the Bach tensor in \cite{caochen:2013} for an $n$-dimensional manifold:
\begin{equation}\label{eq.bachtensorn}
    B_{ij} = \frac{1}{n-3}\nabla^k\nabla^l W_{ikjl} + \frac{1}{n-2} R_{kl}W_{i\, j}^{\, k\, l}
\end{equation}
Clearly, this reduces to the 4-dimensional Bach tensor, although it is not conformally invariant in dimensions $n \geq 5$. We define a gradient shrinking Ricci soliton to be a smooth manifold $(M^n, g)$ which is equipped with a smooth function $f:M^n \rightarrow \mathbb{R}$ satisfying the equation:
\begin{equation}\label{eq.gradientsoliton}
    \mathrm{Hess}(f) + \mathrm{Rc} = \frac{1}{2} g
\end{equation}

H.-D. Cao and Q. Chen use the $D$-tensor defined in \cite{caochen:2013} and \cite{caochen:2009} to show that if $B=0$ on a shrinking or steady soliton, then the $D$-tensor must also vanish. This tensor is given by:
\begin{equation}\label{eq.Dtensor}
    D_{ijk} = \frac{1}{n-2}\left( P_{jk}\nabla_i f - P_{ik}\nabla_j f\right) + \frac{1}{(n-1)(n-2)}\left(g_{jk}E_{il} - g_{ik}E_{jl}\right) \nabla_l f
\end{equation}
H.-D. Cao and Q. Chen then study the implications of the vanishing of the Bach tensor on the geometry of these spaces. Before we can fully motivate our plan in this paper, we would like to express the Bach tensor in terms of the more basic tensors mentioned in the abstract:
    \begin{theorem}[Bergman \cite{bergman:2004}]
        There are only three independent, symmetric, and divergence-free 2-tensors quadratic in the Riemann curvature tensor given by:
        \begin{align}
            \begin{split}\label{eq.utensor}
                U_{ij} &= \, 2(n-3)R_{ipjq}R^{pq} + (n-3)\Delta R_{ij} - \frac{1}{2}(n-3)|Rc|^2 \, g_{ij} \\
                &- (n-3)R \, R_{ij} -\frac{1}{2}(n-3) \Delta R \, g_{ij} + \frac{1}{4}(n-3)R^2 \, g_{ij}
            \end{split}\\
            \begin{split}\label{eq.vtensor}
                V_{ij} =& \,- \nabla_i \nabla_j R + \Delta R \, g_{ij} + R \,
                R_{ij} - \frac{1}{4} R^2 \, g_{ij}
            \end{split}\\
            \begin{split}\label{eq.wtensor}
                W_{ij} =& \, \tensor{R}{_i^{pqr}}\tensor{R}{_{jpqr}} - \frac{1}{4} |Rm|^2 \, g_{ij} - 2R_{ipjq}R^{pq} \\
                &+ R\, R_{ij} - 2R_{pi} R^p_j + |Rc|^2 \, g_{ij} -\frac{1}{4}R^2 \, g_{ij}.
            \end{split}
        \end{align}
    \end{theorem}
    In dimension 4 the $W_{ij}$ tensor vanishes, so we will not consider it here. Also, for clarity of usage we will express $U$ and $V$ in dimension 4 (although the expression for $V$ doesn't change):
    \begin{align}
        \begin{split}\label{eq.utensor4}
            U_{ij} &= \, 2R_{ipjq}R^{pq} + \Delta R_{ij} - \frac{1}{2}|Rc|^2 \, g_{ij} - R \, R_{ij} -\frac{1}{2} \Delta R \, g_{ij} + \frac{1}{4}R^2 \, g_{ij}
        \end{split}\\
        \begin{split}\label{eq.vtensor4}
            V_{ij} &= \,- \nabla_i \nabla_j R + \Delta R \, g_{ij} + R \,
            R_{ij} - \frac{1}{4} R^2 \, g_{ij}
        \end{split}
    \end{align}
    Using the properties of $U$ and $V$ gives us the following expression for the Bach tensor:
    \begin{proposition}
        The Bach tensor given by \eqref{eq.bachtensor} can be expressed by:
        \begin{equation}
            B_{ij} = \frac{1}{2}U_{ij} + \frac{1}{6} V_{ij}
        \end{equation}
        in dimension 4.
    \end{proposition}
    \begin{proof}
        Since $U$ and $V$ are the only two symmetric, quadratic in Riemann curvature, divergence-free tensors in dimension 4, we must have:
        \begin{equation*}
            B_{ij} = \alpha U_{ij} + \beta V_{ij}
        \end{equation*}
        for some constants $\alpha$ and $\beta$. Since the Bach tensor is trace-free, by noticing that the traces of $U$ and $V$ are simply multiples of $\Delta R$, and as $B$ is trace-free:
        \begin{equation*}
            0 = tr(B) = -\alpha \Delta R + 3\beta \Delta R
        \end{equation*}
        So we have $3\beta - \alpha = 0$. Then just by comparing coefficients in the Bach tensor (expressed in terms of $Rc$, $R$, etc., from $W_{ijkl}$), we have $\alpha = \frac{1}{2}$ and $\beta = \frac{1}{6}$.
    \end{proof}
    A question brought up by H.-D. Cao was: if we take another combination of $U$ and $V$, and this tensor vanishes, can we say anything about the geometry of a soliton? The answer turns out to be \textbf{yes}, for some combinations. We will define a Bach-like tensor to be of the following form:
    \begin{equation}\label{eq.bachlike}
            \mathfrak{B}_{ij} = \alpha  U_{ij} + \beta  V_{ij},
        \end{equation}
        where $\alpha$ and $\beta$ are constants.
    Our program for understanding $\mathfrak{B}=0$ on shrinking solitons will first consider the simpler $V$ tensor independently. We construct some integral identities, and then examine the implications of $V=0$. We will then use similar arguments to understand the $U$ tensor indecently.
\section{Properties and Identities for U and V}
    We begin by noticing that the Bach tensor is exactly the right combination of $U$ and $V$ to be trace-free. To see exactly why that is, we simply compute the following trace identities:
    \begin{proposition}\label{prop.VtensorTrace}
        The trace of the tensor $V$ given by \eqref{eq.vtensor} in dimension 4 is $tr(V) = 3\Delta R$, and if $V\equiv 0$ then the scalar curvature is harmonic, i.e. $\Delta R = 0$.
    \end{proposition}
    \begin{proof}
        We simply take the trace of $V$:
        \begin{align*}
            g^{ij}V_{ij} &= -\Delta R + 4\Delta R - R^2 + \frac{1}{4}R^2 (4)
            \\
            &= 3\Delta R.
        \end{align*}
        Then if $V = 0$, we have $\Delta R = 0$.
    \end{proof}
     \begin{proposition}\label{prop.Utensortrace}
        The trace of the tensor $U$ given by \eqref{eq.utensor} in dimension 4 is $tr(U) = -\Delta R$, and if $U\equiv 0$ then the scalar curvature is harmonic, i.e. $\Delta R = 0$.
    \end{proposition}
    \begin{proof}
        We directly compute:
        \begin{equation*}
            g^{ij}U_{ij} = 2|Rc|^2 + \Delta R - 2|Rc|^2 - R^2 - 2\Delta R + R^2.
        \end{equation*}
        So our claim follows.
    \end{proof}
    \begin{lemma}\label{lemma.usefulintegrals}
        Let $(M^4, g, f)$ be a complete gradient shrinking Ricci soliton satisfying \eqref{eq.gradientsoliton} and assume that $\Delta R =0$. Let $r$ be a regular value of $f$. Letting $\Omega = \Omega_r = \{ x \in M \, | \, f(x) \leq r \}$, which is compact as $f$ is proper, we have:
        \begin{enumerate}
            \item[1.)] $\int_{\partial \Omega} \frac{1}{|\nabla f|}\bra \nabla R, \nabla f \ket \diff{S} = 0.$
            
            \item[2.)] $\int_{\Omega}\bra \nabla R, \nabla f \ket \diff{V} = 0.$

            \item[3.)] $\int_{\Omega} |\nabla R|^2 \diff{V} = \int_{\partial \Omega} \frac{R}{|\nabla f|}\bra \nabla R, \nabla f \ket \diff{S} = -\int_{\partial \Omega}|\nabla f| \bra \nabla R, \nabla f \ket \diff{S}.$
            
            \item[4.)] $\int_{\Omega} \bra \nabla R, \nabla f\ket \, e^{-f} \diff{V} = 0.$
            
            \item[5.)] $\int_{\Omega} R Rc(\nabla f, \nabla f)\, e^{-f} \diff{V} = \frac{1}{2}\int_{\partial \Omega}\Big( |\nabla f| \bra \nabla R, \nabla f \ket e^{-f} \diff{S} + \frac{1}{2}\int_{\Omega} |\nabla R|^2 \Big)\, e^{-f} \diff{V}.$
            
            \item[6.)] $\int_{\Omega} f\bra \nabla R, \nabla f \ket \, e^{-f} \diff{V} = 0.$
        \end{enumerate}
    \end{lemma}
     \begin{proof}
        First, we have:
        \begin{equation*}
            0 = \int_{\Omega} \Delta R \diff{V} = \int_{\partial \Omega} \frac{1}{|\nabla f|} \bra \nabla R, \nabla f \ket \diff{S}
        \end{equation*}
        showing (1). Now we integrate by parts to get:
        \begin{align*}
            \int_{\Omega} \bra \nabla R, \nabla f \ket \diff{V} &= \int_{\partial \Omega} \frac{f}{|\nabla f|} \bra \nabla R, \nabla f \ket \diff{S} - \int_{\Omega} f \Delta R \diff{V}
            \\
            &= r\int_{\partial \Omega}\frac{1}{|\nabla f|}\bra \nabla R, \nabla f \ket \diff{S} - \int_{\Omega} f \Delta R \diff{V}.
        \end{align*}
        As $\Delta R = 0$ by assumption, and using (1) we see that both integrals must vanish, showing (2). To show (3) we integrate $R\Delta R$:
        \begin{equation*}
            0 = \int_{\Omega} R\Delta R \diff{V} = \int_{\partial \Omega} \frac{R}{|\nabla f|} \bra \nabla R, \nabla f \ket \diff{S} - \int_{\Omega} |\nabla R|^2 \diff{V}.
        \end{equation*}
        And by using the identity $R = f - |\nabla f|^2$ from normalizing $f$ as in \cite{hamilton:1988}, we get (3). Now we have:
        \begin{align*}
            0 &= \int_{\Omega} \Delta R \, e^{-f} \diff{V}
            \\
            &= \int_{\partial \Omega} \frac{e^{-f}}{|\nabla f|}\bra \nabla R, \nabla f \ket \diff{S} + \int_{\Omega} \bra \nabla R, \nabla f \ket \, e^{-f} \diff{V}
            \\
            &= e^{-r}\int_{\partial \Omega} \frac{1}{|\nabla f|}\bra \nabla R, \nabla f \ket \diff{S} + \int_{\Omega} \bra \nabla R, \nabla f \ket \, e^{-f} \diff{V}
            \\
            &= \int_{\Omega} \bra \nabla R, \nabla f \ket \, e^{-f} \diff{V}
        \end{align*}
        showing (4). To show (5) we integrate $R\Delta R$ with the weighted volume $e^{-f} \diff{V}$:
        \begin{align*}
            0 &= \int_{\Omega} R\Delta R \, e^{-f} \diff{V}
            \\
            &= \int_{\partial\Omega} \frac{R}{|\nabla f|} \bra \nabla R, \nabla f \ket \, e^{-f} \diff{S} - \int_{\Omega} \bra \nabla \big(R e^{-f}\big), \nabla R \ket \diff{V}
            \\
            &= \int_{\partial\Omega} \frac{R}{|\nabla f|} \bra \nabla R, \nabla f \ket \, e^{-f} \diff{S} - \int_{\Omega} \Big( |\nabla R|^2 - R \bra \nabla R, \nabla f \ket \Big) \, e^{-f} \diff{V}
            \\
            &= -\int_{\partial\Omega} |\nabla f| \bra \nabla R, \nabla f \ket \, e^{-f} \diff{S} - \int_{\Omega} \Big( |\nabla R|^2 - 2R Rc(\nabla f, \nabla f) \Big) \, e^{-f} \diff{V}.
        \end{align*}
        And lastly for (6), we have:
        \begin{align*}
            \int_{\Omega} f\bra \nabla R, \nabla f \ket \, e^{-f} \diff{V} &= - \int_{\Omega} \bra f \nabla R, \nabla e^{-f} \ket \diff{V}
            \\
            &= -\int_{\partial \Omega} \frac{f}{|\nabla f|}\bra \nabla R, \nabla f \ket \, e^{-f} \diff{S} + \int_{\Omega} div(f \nabla R) e^{-f} \diff{V}
            \\
            &= -re^{-r}\int_{\partial \Omega} \frac{1}{|\nabla f|}\bra \nabla R, \nabla f \ket  \diff{S} + \int_{\Omega} div(f \nabla R) e^{-f} \diff{V}
            \\
            &= \int_{\Omega} \Big( \bra \nabla R, \nabla f \ket + f \Delta R\Big) e^{-f} \diff{V}
            \\
            &= 0.
        \end{align*}
    \end{proof}
\section{Vanishing V-Tensors on Solitons}
    We can now begin looking at the situation where $V=0$. We have the following:
    \begin{lemma}\label{lemma.intergralV}
        Let $(M^4, g, f)$ be a four dimensional complete gradient shrinking Ricci soliton satisfying \eqref{eq.gradientsoliton} with $V = 0$. Then for any level set $\Omega = \{x\in M^4 \, | \, f(x) \leq r\} \subset M$, we have:
        \begin{equation}\label{eq.vintegral}
            \int_{\Omega} V(\nabla f, \nabla f) \diff{V} = \frac{1}{2} \int_{\Omega}\Big(\nabla R - \frac{R}{2}\nabla f|^2 - \frac{3}{4}R^2 |\nabla f|^2\Big) \diff{V} = 0
        \end{equation}
    \end{lemma}
    \begin{proof}
        We break $\int_{\Omega} V(\nabla f, \nabla f) \diff{V}$ into four parts, evaluate these integrals separately, and then combine them:
        \begin{align*}
           I_1 &= -\int_{\Omega} \nabla^2 R(\nabla f, \nabla f) \diff{V}\\
           I_2 &= \int_{\Omega} |\nabla f|^2 \Delta R \diff{V}\\
           I_3 &= \int_{\Omega} R \, Rc(\nabla f, \nabla f) \diff{V}\\
           I_4 &= -\int_{\Omega} \frac{1}{4} |\nabla f|^2 R^2 \diff{V}.
        \end{align*}
        For the first integral we have:
        \begin{equation*}
            I_1 = -\int_{\Omega}\Big( \bra \nabla \bra \nabla R, \nabla f \ket, \nabla f \ket - Rc(\nabla f, \nabla f) + Rc(\nabla f, \nabla R)\Big) \diff{V}.
        \end{equation*}
        Then by integration by parts and Lemma \ref{lemma.usefulintegrals}(3) we have:
        \begin{align*}
            I_1 &= -\int_{\partial \Omega} |\nabla f| \bra \nabla R, \nabla f \ket \diff{S} + \int_{\Omega}\Big( \bra \nabla R, \nabla f \ket \Delta f + Rc(\nabla f, \nabla f) - Rc(\nabla f, \nabla R) \Big) \diff{V}
            \\
            &=\int_{\Omega}|\nabla R|^2 \diff{V} + \int_{\Omega} \Big( (2 - R)\bra \nabla R, \nabla f \ket  + \frac{1}{2}\bra \nabla R, \nabla f \ket - \frac{1}{2}|\nabla R|^2 \Big)\diff{V}
            \\
            &=\int_{\Omega}\frac{1}{2} \Big(|\nabla R|^2 + \frac{3-2R}{2}\bra \nabla R, \nabla f \ket\Big) \diff{V}
        \end{align*}
        Now $I_2$ simply vanishes under the assumption that $V=0$ by Proposition \ref{prop.VtensorTrace}. The third integral becomes:
        \begin{equation*}
            I_3 = \frac{1}{2}\int_{\Omega} R\bra \nabla R, \nabla f \ket \diff{V}.
        \end{equation*}
        Now combining we have:
        \begin{align*}
            \int_{\Omega} V(\nabla f, \nabla f) \diff{V} &= I_1 + I_2 + I_3 + I_4
            \\
            &= \frac{1}{2}\int_{\Omega}\Big( |\nabla R|^2 + (3 - R)\bra \nabla R, \nabla f \ket - \frac{1}{2}|\nabla f|^2 R^2 \Big)\diff{V}
        \end{align*}
        We note that by Lemma \ref{lemma.usefulintegrals}(2):
        \begin{equation*}
            \int_{\Omega} (3-R)\bra \nabla R, \nabla f \ket \diff{V} = -\int_{\Omega} R\bra \nabla R, \nabla f \ket \diff{V}
        \end{equation*}
        Finally we can complete the square to arrive at:
        \begin{align*}
            \int_{\Omega} V(\nabla f, \nabla f) \diff{V} &= \frac{1}{2}\int_{\Omega}\Big( |\nabla R|^2 - R\bra \nabla R, \nabla f \ket - \frac{1}{2}R^2 |\nabla f|^2 \Big) \diff{V}
            \\
            &= \frac{1}{2} \int_{\Omega}\Big(|\nabla R - \frac{R}{2}\nabla f|^2 - \frac{3}{4}R^2 |\nabla f|^2\Big) \diff{V}.
        \end{align*}
    \end{proof}
    From this it's difficult to say much about the geometry of $M^4$. However, the weighted version of the integral turns out to be more fruitful. For this, we first have to gain control over some terms and do a bit more computing.
    \begin{lemma}\label{lemma.hessianRintegral}
        Let $(M^4, g, f)$ be a complete gradient shrinking Ricci soliton satisfying \eqref{eq.gradientsoliton} and assume that $\Delta R = 0$. Then if $r$ is a regular value of $f$ and $\Omega = \Omega_r = \{ x \in M \, | \, f(x) \leq r \}$, we have:
        \begin{equation}
            \int_{\Omega} \nabla^2 R(\nabla f, \nabla f) \, e^{-f} \diff{V} = -e^{-r} \int_{\Omega} |\nabla R|^2  \diff{V} + \frac{1}{2}\int_{\Omega} |\nabla R|^2 \, e^{-f} \diff{V}.
        \end{equation}
    \end{lemma}
    \begin{proof}
        We integrate:
        \begin{align*}
            \int_{\Omega} \nabla^2 &R(\nabla f, \nabla f) \, e^{-f} \diff{V} 
            \\ &= \int_{\Omega} \bra \nabla \bra \nabla R, \nabla f \ket, \nabla f \ket \, e^{-f} \diff{V} - \frac{1}{2}\cancel{\int_{\Omega} \bra \nabla R, \nabla f \ket \, e^{-f} \diff{V}}\\
            &\hspace{0.5in}+ \frac{1}{2} \int_{\Omega} |\nabla R|^2 \, e^{-f} \diff{V}
            \\
            &= -\int_{\Omega} \bra \nabla \bra \nabla R, \nabla f \ket, \nabla e^{-f} \ket \diff{V} + \frac{1}{2} \int_{\Omega} |\nabla R|^2 \, e^{-f} \diff{V}
            \\
            &= \int_{\partial \Omega} |\nabla f| \bra \nabla R, \nabla f \ket e^{-f} \diff{S} + \int_{\Omega} \Delta (e^{-f}) \bra \nabla R, \nabla f \ket \diff{V}\\
            &\hspace{0.5in}+ \frac{1}{2}\int_{\Omega} |\nabla R|^2 \, e^{-f} \diff{V}
            \\
            &= \int_{\partial \Omega} |\nabla f| \bra \nabla R, \nabla f \ket e^{-f} \diff{S} + \int_{\Omega} (|\nabla f|^2 - \Delta f) \bra \nabla R, \nabla f \ket e^{-f} \diff{V}\\
            &\hspace{0.5in}+ \frac{1}{2}\int_{\Omega} |\nabla R|^2 \, e^{-f} \diff{V}
            \\
            &= \int_{\partial \Omega} |\nabla f| \bra \nabla R, \nabla f \ket e^{-f} \diff{S} + \cancel{\int_{\Omega} (f - 2)\bra \nabla R, \nabla f \ket e^{-f}\diff{V}}\\
            &\hspace{0.5in}+ \frac{1}{2} \int_{\Omega} |\nabla R|^2 \, e^{-f}\diff{V},
        \end{align*}
        giving us the desired form.
    \end{proof}
    \begin{lemma}\label{lemma.gradRbounded}
        Let $(M^n,g_{ij},f)$ be a complete gradient shrinking Ricci soliton satisfying \eqref{eq.gradientsoliton}. Then:
        \begin{equation*}
            \int_M |\nabla R|^2 e^{-\alpha f}\diff{V} < \infty
        \end{equation*}
        for any $\alpha > 0$.
    \end{lemma}
    \begin{proof}
        By the Cauchy-Schwarz inequality we have:
        \begin{equation*}
            \frac{1}{4}|\nabla R|^2 = |Rc(\nabla f)|^2 \leq |Rc|^2 |\nabla f|^2
        \end{equation*}
        By \cite{caozhou:2010}, we know that:
        \begin{equation*}
            |\nabla f|^2(x) \leq \frac{1}{4}(d(x) + c_2)^2
        \end{equation*}
        where $d(x)$ is the distance function from a fixed point $x_0$ and $c_1,c_2$ are constants depending on $n$ and the geometry of $g_{ij}$ of the unit ball $B_{x_0}(1)$. Hence, there exists a compact set $K$ where $|\nabla f|^2 e^{\frac{\alpha f}{2}} \leq 1$ on $M\backslash K$. We have
        \begin{equation*}
            \int_M |Rc|^2 |\nabla f|^2 e^{-\alpha f} \diff{V}= \int_K |Rc|^2 |\nabla f|^2 e^{-\alpha f} \diff{V} +  \int_{M \backslash K} |Rc|^2 |\nabla f|^2 e^{-\alpha f} \diff{V}
        \end{equation*}
        Clearly the first integral over the compact set $K$ is finite. Now we have for the second integral:
        \begin{equation*}
            \begin{split}
                \int_{M \backslash K} |Rc|^2 |\nabla f|^2 e^{-\alpha f} \diff{V} &=  \int_{M \backslash K} (|Rc|^2 e^{-\alpha f/2}) |\nabla f|^2 e^{-\alpha f/2}\diff{V}\\
                &\leq  \int_{M \backslash K} |Rc|^2 e^{-\alpha f/2}\diff{V}
            \end{split}
        \end{equation*}
        By O. Munteanu and N. Sesum's results in \cite{munteanusesum:2013}, we have:
        \begin{equation*}
            \int_M |Rc|^2 e^{-\lambda f} \diff{V} < \infty
        \end{equation*}
        for any $\lambda > 0$, so after setting $\lambda = \frac{\alpha}{2}$ the result follows.
    \end{proof}
    \begin{lemma}\label{lemma.gradRvanishes}
        Let $(M^n,g_{ij},f)$ be a complete gradient shrinking Ricci soliton satisfying \eqref{eq.gradientsoliton}. Then if $\Omega_r = \{ f\leq r\}$  for regular values $r > 0$ of $f$ and $\alpha >0$:
        \begin{equation*}
            \lim_{r\rightarrow \infty} \Big(e^{-\alpha r}\int_{\Omega_r} |\nabla R|^2 \diff{V}\Big) = 0
        \end{equation*}
    \end{lemma}
    \begin{proof}
    By Lemma \ref{lemma.gradRbounded} we know that:
    \begin{equation*}
        \int_M |\nabla R|^2 e^{-\alpha f/2} \diff{V} < \infty
    \end{equation*}
    Since $f \leq r$ on $\Omega_r$, we also have $e^{-\alpha r/2} \leq e^{-\alpha f/2}$, giving us:
    \begin{equation*}
        \frac{e^{-\alpha r}}{2}\int_{\Omega_r} |\nabla R|^2 \diff{V} = \frac{e^{-\alpha r/2}}{2} \int_{\Omega_r} |\nabla R|^2 e^{-\alpha r/2}\diff{V} \leq \frac{e^{-\alpha r/2}}{2}\int_{\Omega_r} |\nabla R|^2 e^{-\alpha f/2}\diff{V}
    \end{equation*}
    Letting $r\rightarrow \infty$, we see that
    \begin{equation*}
        \lim_{r\rightarrow \infty}\Big(\frac{e^{-\alpha r}}{2}\int_{\Omega_r}|\nabla R|^2 \diff{V}\Big) \leq \lim_{r\rightarrow \infty}\Big(\frac{e^{-\alpha r/2}}{2}\int_{\Omega_r} |\nabla R|^2 e^{-\alpha f/2}\diff{V}\Big) = 0
    \end{equation*}
    \end{proof}
    We are now ready to prove a useful theorem about the geometry of a complete gradient shrinking Ricci soliton if the $V$-tensor vanishes.
    \begin{theorem}\label{theorem.vanishingV}
        Let $(M^4, g, f)$ be a complete gradient shrinking Ricci soliton satisfying \eqref{eq.gradientsoliton}. Then $V\equiv 0$ if and only if $M^4$ is isometric to an Einstein manifold or a Gaussian soliton.
    \end{theorem}
    \begin{proof}
    Let $\Omega = \Omega_r = \{ x \in M \, | \, f(x) \leq r \}$ where $r$ is a regular value of $f$. By Proposition \ref{prop.VtensorTrace}, $V\equiv 0$ implies $\Delta R = 0$. By Lemma \ref{lemma.usefulintegrals}(5) and (3) we have:
    \begin{align*}
         \int_{\Omega} R\, Rc(\nabla f, \nabla f) \, e^{-f} \diff{V} &= -\frac{1}{2}\int_{\Omega} R\Delta_f R e^{-f}\diff{V}
        \\
        &= -\frac{1}{2}\int_{ \Omega} |\nabla R|^2 \, e^{-r} \diff{V} + \frac{1}{2}\int_{\Omega} |\nabla R|^2 \, e^{-f} \diff{V}.
    \end{align*}
    And using Lemma \ref{lemma.hessianRintegral}, we compute the integral of $V(\nabla f, \nabla f)$ with respect to the weighted volume element $e^{-f}\diff{V}$, and we note that the second term involving $\Delta R$ again vanishes:
    \begin{align*}
        \int_{\Omega} V(&\nabla f, \nabla f) \, e^{-f} \diff{V} \\
        &= -\int_{\Omega} \Big(\nabla^2 R(\nabla f, \nabla f) - R \, Rc(\nabla f, \nabla f) + \frac{1}{4}R^2 |\nabla f|^2 \Big)\, e^{-f} \diff{V}
        \\
        &= \frac{e^{-r}}{2}\int_{\Omega} |\nabla R|^2 \diff{V} - \frac{1}{4}\int_{\Omega}R^2|\nabla f|^2 \, e^{-f} \diff{V}.
    \end{align*}
    Letting $r\rightarrow \infty$, by Lemma \ref{lemma.gradRvanishes} the first integral vanishes, so we see that:
    \begin{equation}
        \label{eq.VweightedIntegral}
        \int_M V(\nabla f, \nabla f)\, e^{-f}\diff{V} = -\int_{M}\frac{1}{4}R^2 |\nabla f|^2 \, e^{-f} \diff{V} 
    \end{equation}
    so that if $V=0$, this integral must vanish.
    Therefore at any point $p\in M$ we must have either $R = 0$ or $|\nabla f| = 0$. By the work of S. Pigola, M. Rimoldi, and A. Setti \cite{pigolarimoldi:2011}, and say either $R>0$ everywhere, hence $\nabla f = 0$ and the manifold is Einstein, or if $R=0$ anywhere, then the soliton must be isometric to the Gaussian soliton on $\mathbb{R}^4$. But we can see most of this result directly: using Lemma \ref{lemma.intergralV} we must have:
    \begin{equation*}
        \frac{1}{2} \int_{M}\Big(|\nabla R - \cancel{\frac{R}{2}\nabla f}|^2 - \cancel{\frac{3}{4}R^2 |\nabla f|^2}\Big) \diff{V} = \frac{1}{2}\int_{M}|\nabla R|^2 = 0.
    \end{equation*}
    Hence, the scalar curvature is constant. Moreover, if $R\neq 0$, then $\nabla f$ must vanish on all of $M$, showing the manifold is Einstein. In fact  we have $R=2$ since:
    \begin{equation*}
        R = \frac{1}{2} \mbox{tr}(g) = 2.
    \end{equation*}
    Now if $R=0$, then from the Bochner formula we have:
    \begin{align*}
        \frac{1}{2}\Delta |\nabla f|^2&= \bra \nabla \Delta f, \nabla f \ket + |\nabla^2 f|^2 + Rc(\nabla f, \nabla f)
        \\
        &= \bra \nabla(2-R), \nabla f \ket + |\nabla^2 f|^2 + \frac{1}{2}\bra \nabla R, \nabla f \ket 
    \end{align*}
    Hence, we get:
    \begin{align*}
        1 - \cancel{\frac{R}{2}} &= -\cancel{\frac{1}{2}\bra \nabla R, \nabla f \ket} + |\nabla^2 f|^2
        \\
        1 &= |\frac{1}{2} g - Rc|^2
        \\
        1 &= 1 - R + |Rc|^2
        \\
        0 &= |Rc|^2
    \end{align*}
    So the manifold is Ricci flat, hence Einstein. For the converse, suppose $(M^4, g, f)$ is Einstein.
    Then as the scalar curvature is constant on any Einstein manifold, we have:
        \begin{align*}
            V_{ij} &= \frac{1}{4} R^2 g_{ij} - R \, R_{ij}
            \\
            &= -R (R_{ij} - \frac{1}{4} R \, g_{ij})
        \end{align*}
        But under our assumption in dimension 4, we have $R_{ij} = \frac{R}{4} g_{ij}$, hence $V = 0$. If the manifold is Gaussian, then $R=0$ and again $V=0$.
    \end{proof}
\section{Analyzing the U-Tensor}
We would like to use some of the tools developed in the previous sections to analyze the situation where $U= 0$. It turns out that while we obtain some suggestive looking $L_2$-norms, we can't quite prove a similar result to when $V=0$.
\begin{lemma}\label{lemma.Uintegral}
        Let $(M^4, g, f)$ be a complete gradient shrinking Ricci soliton satisfying \eqref{eq.gradientsoliton} normalized so that $R + |\nabla f|^2 = f$. Then if $\Delta R = 0$ we have:
        \begin{equation*}
            \int_M U(\nabla f, \nabla f)\diff{V} = \frac{1}{4}\int_{\Omega}|\nabla f|^2 ( R^2 - 2|Rc|^2) \diff{V}.
        \end{equation*}
    \end{lemma}
    
    \begin{proof}
        For $\Omega = \Omega_r = \{x \in M^4 \, |\, f(x) \leq r \}$, we have:
        \begin{align*}
            \int_{\Omega} \nabla_k& R_{ij} \nabla_k f \nabla_i f \nabla_j f \diff{V} 
            \\
            &= \int_{\Omega} \nabla_k(R_{ij}\nabla_i f\nabla_j f\nabla_k f) \diff{V} - 2\int_{\Omega} R_{ij}( \nabla_i\nabla_k f(\nabla_k f))\nabla_j f \diff{V}\\
            & \hspace{0.15in} - \int_{\Omega} \Delta f \, R_{ij}\nabla_i f\nabla_j f \diff{V}
            \\
            &= \frac{1}{2} \int_{\partial \Omega} |\nabla f| \bra \nabla R, \nabla f \ket \diff{S} - 2 \int_{\Omega} Rc(\nabla f, \frac{1}{2}\nabla f - Rc(\nabla f)) \diff{V}\\
            & \hspace{0.15in} + \int_{\Omega} R \, Rc(\nabla f, \nabla f) \diff{V}
            \\
            &= \int_{\Omega} R\, Rc(\nabla f, \nabla f) \diff{V}.
        \end{align*}
        Hence, as the first term $Rc(\nabla f,\nabla f)$ vanishes when we integrate, the fifth term $\Delta R |\nabla f|^2 = 0$, and keeping in mind the fourth term is just $-\int_{\Omega}R\, Rc(\nabla f, \nabla f)\diff{V}$, we get:
        \begin{equation*}
            \int_{\Omega} U(\nabla f, \nabla f)\diff{V} = \frac{1}{2}\int_{\Omega}\frac{R^2}{2}|\nabla f|^2 - |Rc|^2|\nabla f|^2\diff{V}.
        \end{equation*}
    \end{proof}
    \begin{lemma}\label{lemma.weightedUintegral}
        Let $(M^4, g, f)$ be a complete gradient shrinking Ricci soliton satisfying \eqref{eq.gradientsoliton} normalized so that $R + |\nabla f|^2 = f$. Then if $\Delta R =0$ we have:
        \begin{equation*}
             \int_M U(\nabla f, \nabla f) e^{-f} \diff{V} = \frac{1}{4}\int_{M} |\nabla f|^2 (R^2 - 2|Rc|^2)\, e^{-f} \diff{V}.
        \end{equation*}
    \end{lemma}
    \begin{proof}
         For $\Omega = \Omega_r = \{x \in M^4 \, |\, f(x) \leq r \}$, we have:
       \begin{align*}
            \int_{\Omega} &\nabla_k R_{ij} \nabla_k f \nabla_i f \nabla_j f \, e^{-f} \diff{V}
            \\
            &= \int_{\Omega} \nabla_k (R_{ij} \nabla_i f\nabla_j f\nabla_k f \, e^{-f})\diff{V} - 2\int_{\Omega} Rc(\nabla f, \frac{1}{2}\nabla f - Rc(\nabla f))\, e^{-f} \diff{V}\\
            &\hspace{0.15in} - \int_{\Omega} \Delta f\, Rc(\nabla f, \nabla f) e^{-f}\diff{V} + \int_{\Omega}|\nabla f|^2 Rc(\nabla f, \nabla f) \, e^{-f} \diff{V}
            \\
            &= \int_{\Omega} R \, Rc(\nabla f, \nabla f) \, e^{-f} \diff{V}.
        \end{align*}
        So, after adding everything together we get:
        \begin{equation*}
            \int_{\Omega_r} U(\nabla f, \nabla f) \, e^{-f} \diff{V} =  \frac{1}{4}\int_{\Omega} |\nabla f|^2 (R^2 - 2|Rc|^2)\, e^{-f} \diff{V}.'
        \end{equation*}
    \end{proof}
    It's interesting that the integrand of both the weighted and unweighted integrals are exactly the same. However, the integrand $R^2-2|Rc|^2$ may change sign from point to point.
    \begin{align*}
        \nabla_k (&R_{ij} \nabla_i f\nabla_j f \nabla_k f)\\
        &=\mathrm{div}(\frac{1}{2}\bra \nabla R, \nabla f \ket \nabla f)\\
        &=\frac{1}{2}\bra \nabla\bra \nabla R, \nabla f \ket, \nabla f \ket + \frac{1}{2}\bra \nabla R, \nabla f\ket \Delta f\\\
        &=\frac{1}{2}\hess R (\nabla f, \nabla f) + \frac{1}{2}\hess f(\nabla R, \nabla f) + \bra \nabla R, \nabla f \ket - R \, \Rc(\nabla f, \nabla f)\\
        &= \frac{1}{2}\hess R(\nabla f, \nabla f) + \frac{5}{4}\bra \nabla f, \nabla R \ket - \frac{1}{4}|\nabla R|^2 - R\, \Rc(\nabla f, \nabla f)
    \end{align*}
    \begin{align*}
        \nabla_k (&R_{ij} \nabla_i f\nabla_j f \nabla_k f) 
        \\
        &=(\nabla_k R_{ij}) \nabla_i f\nabla_j f \nabla_k f + 2R_{ij} \nabla_i \nabla_k f\nabla_k f \nabla_j f + (2-R)R_{ij} \nabla_i f\nabla_j f
        \\
        &= (\nabla_k R_{ij}) \nabla_i f\nabla_j f  \nabla_k f + 3\Rc(\nabla f, \nabla f) -\frac{1}{2} |\nabla R|^2- R\,\Rc(\nabla f, \nabla f)
    \end{align*}
    Therefore, we have:
    \begin{equation}
        \begin{split}
            (\nabla_k R_{ij}& \nabla_k) \nabla_i f\nabla_j f\\
            &= \frac{1}{2}\nabla^2 R(\nabla f, \nabla f) - \frac{1}{2}\Rc(\nabla f, \nabla f) +\frac{1}{4}|\nabla R|^2
        \end{split}
    \end{equation}
    We have as an elementary fact that:
    \begin{equation*}
        R^2 - 4|Rc|^2 \leq 0
    \end{equation*}
    Moreover, we know that $R^2 = 4|Rc|^2$ is equivalent to $M^4$ being Einstein. So far, we can only say the following:
    \begin{proposition}
        Let $(M^4,g_{ij},f)$ be an Einstein manifold. Then $U=0$.
    \end{proposition}
    \begin{proof}
        If $M^4$ is Einstein then $Rc = \frac{1}{4}R \, g$, so:
        \begin{align*}
            \frac{1}{n-3}U_{ij} &= 2R_{ipjq}\frac{R}{4}g^{pq} - \frac{1}{2}\Big(\frac{R^2}{4}\Big) g_{ij} - \frac{R^2}{4}g_{ij} + \frac{R^2}{4}g_{ij}\\
            &= \frac{R}{2}R_{ij} - \frac{R^2}{8}g_{ij}
            \\
            &= \frac{R^2}{8}g_{ij} - \frac{R^2}{8}g_{ij} = 0.
        \end{align*}
    \end{proof}
    But as was illustrated above, the converse is less clear. There may be another approach or perhaps we can use some of the techniques developed here in the future to prove the following conjecture (or something similar):
    \begin{conjecture}
        Let $(M^4, g, f)$ be a complete gradient shrinking Ricci soliton satisfying \eqref{eq.gradientsoliton} normalized so that $R + |\nabla f|^2 = f$. Then $U=0$ if and only if $M^4$ is Einstein or isometric to a Gaussian soliton or isometric to $\mathbb{S}^3 \times \mathbb{R}$.
    \end{conjecture}
\section{Vanishing Bach-like Tensors}
    To understand the question posed by H.-D. Cao previously, we will study Bach-like tensors. We can express $\mathfrak{B}$ in terms of the Bach tensor in the following way:
    \begin{equation}\label{eq.bachlikenormalized}
        \mathfrak{B}_{ij} = 2\alpha \, B_{ij} + \Big(\beta - \frac{\alpha}{3}\Big) V_{ij}
    \end{equation}
    If $\alpha = 0$ then $\mathfrak{B} = \beta V$, and we already have results for this situation. We also know the following:
    \begin{lemma}
        Let $(M^4,g,f)$ be a complete gradient shrinking Ricci soliton satisfying \eqref{eq.gradientsoliton}. Then if $\mathfrak{B} = 0$ and $\beta - \alpha/3 \neq 0$, the scalar curvature is harmonic. ie $\Delta R = 0$.
    \end{lemma}
    \begin{proof}
        By Propositions \ref{prop.Utensortrace} and \ref{prop.VtensorTrace}, we have:
        \begin{equation*}
            0 = g^{ij}\Big(\alpha U_{ij} + \beta V_{ij}\Big) = \Big(3\beta - \alpha)\Delta R.
        \end{equation*}
    \end{proof}
    \begin{lemma}\label{lemma.bachLikeFlatIntegral1}
        Let $(M^4,g,f)$ be a complete gradient shrinking Ricci soliton satisfying \eqref{eq.gradientsoliton}. If $\mathfrak{B} = 0$ and $3\beta - \alpha \neq 0$ then 
        \begin{equation*}
            \int_M \Big((\alpha - \beta)|\nabla f|^2 R^2 - 2\alpha|Rc|^2\Big) \, e^{-f}\diff{V} = 0
        \end{equation*}
    \end{lemma}
    \begin{proof}
        Given our assumptions about $\alpha$ and $\beta$, clearly $\mathfrak{B} = 0$ implies that $\Delta R = 0$. Then using \ref{lemma.weightedUintegral} and  \eqref{eq.vintegral}:
        \begin{equation*}
            \alpha\int_M |\nabla f|^2(R^2 - 2|Rc|^2)e^{-f}\diff{V} - \beta\int_M R^2|\nabla f|^2 e^{-f} \diff{V} = 0
        \end{equation*}
        Hence
        \begin{equation*}
            \int_M \Big((\alpha - \beta)|\nabla f|^2 R^2 - 2\alpha|Rc|^2\Big) \, e^{-f}\diff{V} = 0
        \end{equation*}
    \end{proof}
    From this we have the following:
    \begin{proposition}\label{prop.bachLikeClassification1}
        Let $(M^4,g,f)$ be a complete gradient shrinking Ricci soliton satisfying \eqref{eq.gradientsoliton} with $\mathfrak{B} = \alpha U + \beta V = 0$. If $0 \leq \alpha \leq \beta$, or $\alpha \geq \beta$  and $\alpha < 0$, then $M^4$ must be isometric to an Einstein manifold or a Gaussian soliton.
    \end{proposition}
    \begin{proof}
        Simply looking at the second integral from Lemma \ref{lemma.bachLikeFlatIntegral1}, we can use the same arguments from Theorem \ref{theorem.vanishingV}.
    \end{proof}
    So for at least a quarter of the $(U,V)$-plane that Bach-like tensors exist on, if $\mathfrak{B} = 0$ we can classify a complete gradient shrinking Ricci soliton.
    \begin{lemma}\label{lemma.weightedBachIntegral}
        Let $(M^n,g_{ij},f)$ be a complete gradient shrinking Ricci soliton satisfying \eqref{eq.gradientsoliton}. Then
        \begin{equation*}
            \int_M B(\nabla f, \nabla f) e^{-f} \diff{V} = -\frac{1}{2}\int_M |D_{ijk}|^2 e^{-f} \diff{V}.
        \end{equation*}
    \end{lemma}
    \begin{proof}
        Let $r$ be a regular value of $f$ and $\Omega=\Omega_r=\{f\leq r\}$. We use the antisymmetric properties of the Cotton and $D$ tensor:
        \begin{align*}
            (n-2)\int_{\Omega_r} &B_{ij} \nabla_i f\nabla_j f e^{-f}\diff{V}
            \\
            &= -\int_{\Omega_r} \nabla_k D_{ikj} \nabla_i f\nabla_j f e^{-f}\diff{V} - \frac{n-3}{n-2} \int_{\Omega_r} C_{jki} \nabla_i f \nabla_j f \nabla_k  f e^{-f}\diff{V}\\
            &= - \int_{\partial \Omega_r} D_{ikj} \nabla_i f\nabla_j f\nabla_k f \cdot |\nabla f|^{-1}\, e^{-f} \diff{S}  \\
            &\hspace{0.5cm} + \int_{\Omega_r} D_{ikj}(R_{ik}\nabla_j f + R_{kj}\nabla_i f) \, e^{-f}\diff{V}\\
            &\hspace{0.5cm}- \int_M D_{ikj} \nabla_i f\nabla_j f\nabla_k f \, e^{-f}\diff{V}
            \\
            &= \int_{\Omega_r} D_{ikj} (R_{ik} \nabla_j f + R_{kj}\nabla_i f)\, e^{-f} \diff{V}\\
            &= -\frac{1}{2}\int_{\Omega_r} D_{ijk} ( R_{ik}\nabla_j f - R_{kj} \nabla_i f) \, e^{-f} \diff{V}\\
            &= - \frac{(n-2)}{2}\int_{\Omega_r} |D_{ijk}|^2 e^{-f} \diff{V}.
        \end{align*}
    \end{proof}
    \begin{lemma}\label{lemma.bachLikeFlatIntegral2}
        Let $(M^4,g,f)$ be a complete gradient shrinking Ricci soliton satisfying \eqref{eq.gradientsoliton} with $\mathfrak{B} = \alpha U + \beta V = 0$. Then
        \begin{equation*}
             -\alpha \int_M |D_{ijk}|^2 \, e^{-f} \diff{V} - \frac{1}{2}\Big(\beta-\frac{\alpha}{3}\Big)\int_M R^2|\nabla f|^2 \, e^{-f} \diff{V} = 0.
        \end{equation*}
    \end{lemma}
    \begin{proof}
        We have from \eqref{eq.bachlikenormalized} that:
        \begin{equation*}
            \begin{split}
                \int_M \mathfrak{B}&(\nabla f, \nabla f) e^{-f} \diff{V}\\
                &= 2\alpha \int_M B(\nabla f, \nabla f)\, e^{-f} \diff{V} + \Big(\beta - \frac{\alpha}{3}\Big) \int_M V(\nabla f, \nabla f) \, e^{-f} \diff{V}.
            \end{split}
        \end{equation*}
        Now we just replace these integrals using Lemma \ref{lemma.weightedBachIntegral} and \eqref{eq.VweightedIntegral} to complete the proof.
    \end{proof}
    \begin{proposition}\label{prop.bachLikeClassification2}
        Let $(M^4,g,f)$ be a complete gradient shrinking Ricci soliton satisfying \eqref{eq.gradientsoliton} with $\mathfrak{B} = \alpha U + \beta V = 0$. If $(\alpha, \beta) \in \{ (x,y) \, | \, x>0,\, y >x/3\}$ or $(\alpha, \beta) \in \{ (x,y) \, | \, x<0,\, y <x/3\}$ then $M^4$ must be isometric to an Einstein manifold or a Gaussian soliton.
    \end{proposition}
    \begin{proof}
        We simply use Lemma \ref{lemma.bachLikeFlatIntegral2}. If $\alpha > 0$ and $\beta > \alpha / 3$ then it forces $|D_{ijk}|^2 = 0$ and $R^2|\nabla f|^2 = 0$. By the same arguments in Theorem \ref{theorem.vanishingV}, $M^4$ must be Einstein or Gaussian. The situation where $\alpha < 0$ and $\beta < \alpha /3$ is analogous.
    \end{proof}
    Collecting all of this information, as well as what we know from H.-D. Cao and Q. Chen's work on complete gradient shrinking Ricci solitons discussed in the previous chapter and in \cite{caochen:2013}, we have the following:
    \begin{theorem}\label{theorem.bachLikeClassification}
     Let $(M^4,g,f)$ be a complete gradient shrinking Ricci soliton satisfying \eqref{eq.gradientsoliton}, and $\mathfrak{B} = \alpha U + \beta V$ as in \eqref{eq.bachlike}. Let $\Lambda$ be the set given by:
     \begin{equation*}
         \Lambda = \{ (x,y) \, | \, x\geq 0,\, y > x/3\} \cup \{ (x,y) \, | \, x\leq0,\, y < x/3\}
     \end{equation*}
     Assume $\mathfrak{B} = 0$. Then
     \begin{enumerate}
         \item if $(\alpha, \beta) \in \Lambda$, $M^4$ is either Einstein or a Gaussian soliton; or
         \item if $\beta = \alpha/3$, and $\alpha\neq 0$, then $M^4$ is either Einstein, or locally conformally flat, hence a finite quotient of either the Gaussian shrinking soliton $\mathbb{R}^4$ or the round cylinder $\mathbb{S}^3\times \mathbb{R}$.
     \end{enumerate}
    \end{theorem}
    \begin{proof}
        The first claim is just a result from Propositions \ref{prop.bachLikeClassification2}. For the second, we recall that if $\beta = \alpha / 3$ and $\alpha \neq 0$, then $\mathfrak{B} = 2\alpha B$, so the results of H.-D. Cao and Q. Chen in \cite{caochen:2013} hold.
    \end{proof}
    \newpage
    We can visualize the situation in the following way. The tensors $U$ and $V$ are algebraically independent, and the Bach tensor and Bach-like tensors are linear combinations of these. If we imagine the Bach tensor and Bach-like tensors living on the $(U,V)$-plane, we have the following picture:
    
    \vspace{0cm}
    
    \begin{figure}[h]
        \includegraphics[scale=2.5, center]{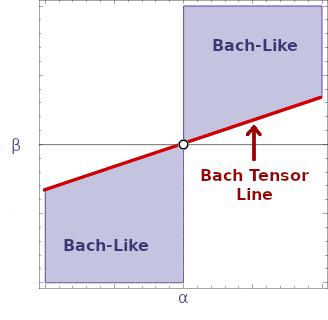}
        \caption{}
        \label{figure.ab-plane}
    \end{figure}
    
    \vspace{0cm}
    
    For combinations of $\alpha$ and $\beta$ along the so-called Bach tensor line in figure \ref{figure.ab-plane}, we have more possibilities, while within the Bach-Like regions, including the case where $\alpha = 0$, we know that $M^4$ is either just Einstein or Gaussian.
\section{Current Work}

\printbibliography
\end{document}